\documentclass[11pt,a4paper]{article}

\usepackage[utf8]{inputenc}
\usepackage[T1]{fontenc}
\usepackage[margin=2.5cm]{geometry}
\usepackage{amsmath,amssymb,amsthm}
\usepackage{graphicx}
\usepackage{booktabs}
\usepackage{multirow}
\usepackage{xcolor}
\usepackage{enumitem}
\usepackage{hyperref}
\usepackage{bbm}

\newcommand{\R}{\mathbb{R}}
\newcommand{\N}{\mathbb{N}}
\newcommand{\cK}{\mathcal{K}}

\newtheorem{theorem}{Theorem}
\newtheorem{proposition}{Proposition}

\theoremstyle{definition}

\title{Generating health variables for a synthetic population of Santa Catarina, Brazil, via ensemble learning and a calibrated linear program}

\author{%
Maria Eduarda Pinheiro\thanks{Departamento de Matem\'atica, Universidade Federal de Santa Catarina, Blumenau, SC, Brazil. \\ maria.eduarda.pinheiro@ufsc.br, vitoria.blau@grad.ufsc.br}
\and
Vitoria Caroline Blau\footnotemark[1]
}
\date{}

\begin{document}
\maketitle

\begin{abstract}
Agent-based models in epidemiology require synthetic populations with
health-related attributes, which in Brazil are available only from sample
surveys such as the Pesquisa Nacional de Sa\'ude (PNS). We add health
variables, including smoking, alcohol consumption, physical exercise,
self-rated health and chronic conditions, to a georeferenced synthetic
population of the state of Santa Catarina. Class probabilities are
estimated by an ensemble of logistic regression, random forest and support
vector machine, and individuals are assigned to classes by a linear program
that reproduces the class proportions estimated from the survey. We show
that, under a mild condition satisfied by our ensemble, the penalty term of
this program is no longer needed, and the proportions are reproduced
exactly. In a cross-validation on PNS data, the method attains the best results.
\end{abstract}

\section{Introduction}
\label{sec:intro}

Synthetic populations are the basic input of agent-based models in
epidemiology, transportation and urban planning
\cite{chapuis2022review,eubank2004}. For epidemiological modeling, their
usefulness depends on attributes that go beyond demographics: smoking,
alcohol consumption and physical inactivity are risk factors for several
chronic diseases, and these factors, together with general health status
and existing chronic conditions, modulate the individual response to
infections and interventions.

For the state of Santa Catarina (SC), Brazil, a georeferenced synthetic
population was recently constructed by combining 2022 Census aggregates
with the 2010 Census microdata through Iterative Proportional Fitting and
integer programming, and then adjusted to 2025 municipal population
estimates \cite{pinheiro2026synthpop}. It contains about 8.2 million
individuals in 3.0 million households, each located at a residence, but
only demographic attributes: age, sex, race/skin color and household
composition. In Brazil, individual-level health information is collected by
the PNS \cite{ibge2020pns,stopa2020pns}, a household survey whose sample
for a single state has only a few thousand respondents.

Transferring survey information to a full population is a classification
problem with two specific difficulties. First, the outcomes are
categorical, often multiclass, and strongly imbalanced: in the SC sample,
for instance, only 2.1\% of respondents report a previous stroke. Second,
what matters for a synthetic population is not only individual-level
accuracy but also that the class proportions of each outcome reproduce the
survey estimates. A classifier that assigns each individual to its most
probable class tends to under-represent rare classes and offers no control
over these proportions.

Burgard et al.\ \cite{burgard2026fair} addressed this problem with a linear
program that assigns classes to unlabeled individuals by balancing their
estimated class probabilities against prescribed class counts, penalizing
deviations from these counts. They proved that the constraint matrix is
totally unimodular, so that the linear relaxation yields an integral
assignment. In this paper, we apply this approach to Brazilian data in two
stages: class probabilities are first estimated by an ensemble of
multinomial logistic regression, random forest and support vector machine,
and individuals are then assigned to classes by the linear program, solved
separately for each municipality.

Our contributions are as follows. First, we show that the ensemble
probabilities are strictly positive, a property inherited from the logistic
regression, and that this positivity, together with target counts that sum
to the number of individuals, implies that every optimal solution of the
linear program attains the target counts exactly whenever the penalty
parameter exceeds $1/m$. The program therefore reduces to a transportation
problem, and the penalty parameter plays no role in the solution. Second,
we evaluate the complete procedure by repeated cross-validation on PNS
data, comparing it with the individual classifiers and with the ensemble
without the assignment step. In this evaluation, the target counts are
taken from the held-out labels, which isolates the effect of the assignment
step. Third, we apply the procedure to generate health variables for the
synthetic population of Santa Catarina and make them publicly available.

The paper is organized as follows. Section~\ref{sec:data} describes the PNS
data, the synthetic population and the generation procedure.
Section~\ref{sec:methods} presents the classifiers, the ensemble, the
optimization model and the structure of its optimal solutions.
Section~\ref{sec:experiments} describes the experimental design,
Section~\ref{sec:results} reports the validation results,
Section~\ref{sec:generation} the generation for the synthetic population,
and Section~\ref{sec:conclusion} concludes.

\section{Data}
\label{sec:data}

\subsection{The National Health Survey (PNS 2019)}

The PNS is a household survey conducted by the Brazilian Institute of
Geography and Statistics (IBGE) in partnership with the Ministry of Health
\cite{ibge2020pns,stopa2020pns}. It has a complex sampling design
(stratified, multi-stage cluster sampling) and provides sampling weights
calibrated to population totals. In Santa Catarina, the survey covers about
$10\,000$ residents. 
The individual questionnaire, which contains the variables used here, is
answered by one resident aged 15 or older, selected at random with equal
probability among the eligible residents of each household
\cite{stopa2020pns}.
After removing records with missing values in the predictors (age, sex and
race) or in the lifestyle and self-rated health variables, the working sample
contains $N=3\,738$ respondents. For diabetes and high cholesterol, records
with missing or inapplicable answers to the corresponding item are also
removed, which leaves $3\,539$ and $3\,514$ respondents, respectively.

\subsection{Target outcomes}
We consider smoking, alcohol consumption, and physical exercises, which also serve as predictors for the self-rated-health outcome. And we consider self-reported diagnoses of chronic diaseases: diabetes, high cholesterol, heart disease, stroke, asthma, arthritis, depression and chronic kidney disease.
Table~\ref{tab:targets} lists information about these outcomes.

\begin{table}[h]
\centering
\caption{Target outcomes from the PNS 2019, SC subsample.} 
\label{tab:targets}
\begin{tabular}{@{}lll@{}}
\toprule
Outcome (PNS item) & Classes & Proportion (\%) \\
\midrule
Smoking (P050) & daily / occasional / does not smoke & 12.5 / 1.1 / 86.4 \\
Alcohol (P027) & never / $<$1 per month / $\geq$1 per month & 52.1 / 13.5 / 34.5 \\
Physical exercise (P034) & yes / no (last 3 months) & 44.4 / 55.6 \\
Self-rated health (N001) & very good / good / fair / poor / very poor & 19.8 / 53.3 / 22.1 / 3.9 / 0.9 \\
Diabetes (Q03001) & yes/ no & 8.8 / 91.2 \\
High cholesterol (Q060) & yes / no & 16.1 / 83.9 \\
Heart disease (Q06306) & yes / no &  7.1 / 92.9 \\
Stroke (Q068) & yes / no &  2.1 / 97.9 \\
Asthma (Q074)                  & yes / no &  5.7 / 94.3 \\
Arthritis/rheumatism (Q079)    & yes / no &  7.9 / 92.1 \\
Depression (Q092) & yes / no & 13.4 / 86.6 \\
Chronic kidney disease (Q124)  & yes / no &  2.1 / 97.9 \\
\bottomrule
\end{tabular}
\end{table}

\subsection{Synthetic population of Santa Catarina}

The synthetic population of \cite{pinheiro2026synthpop} contains
8,187,029 individuals in 3,014,290 households across 295
municipalities, each with age, sex, race,
household identifier, intra-household order and georeferenced residence.

It was built by combining 2022 Census aggregates (household and population totals per census tract) with the household co-residence structure observed in the 2010 Census microdata sample, reconciled through Iterative Proportional Fitting and integer optimization; the resulting 2022 population was then adjusted to match 2025 municipal population estimates through targeted removal and cloning of households in municipalities with, respectively, declining and growing population.
\subsection{Generation procedure}
\label{sec:procedure}

We add to the synthetic population of Santa Catarina generated by \cite{pinheiro2026synthpop} the
health variables of Table~\ref{tab:targets}. Since the PNS individual questionnaire is answered by
residents aged 15 or older, the variables are generated for all heads of household and for the other
household members aged 15 or older.

\paragraph{Predictors.} Only attributes present in both data sources can
be used as predictors. At the individual level these are age, sex and
race. Variables generated in earlier steps are added as
predictors in later steps, as described below. For household members other than the head, the value of the same variable
generated for the head of their household is also used as a predictor.

\paragraph{Target counts.}
Given a variable with $\mathcal K:=\{1,\dots,K\}$ classes, let $\hat\pi_k$ be
the proportion of class $k$ among the PNS respondents of Santa Catarina who are
heads of household, without sampling weights. For a municipality with
$m=N^{(1)}$ heads of household, the cardinality of each class $k$ among the
heads is set to
\begin{equation}
\lambda^{(1)}_k=\operatorname{round}\left(\hat\pi_k\,N^{(1)}\right),
\qquad k\in\mathcal K,
\label{eq:lambda}
\end{equation}
and, for the $m=N^{(j)}$ household members of order $j\ge2$ in the same
municipality,
\begin{equation}
\lambda^{(j)}_k=\operatorname{round}\left(\frac{N^{(j)}}{N^{(1)}}\,
\lambda^{(1)}_k\right), \qquad k\in\mathcal K.
\label{eq:lambdaj}
\end{equation}
In both cases the rounding is done by the largest-remainder method
(Section~\ref{sec:generation}), so that $\sum_k\lambda^{(j)}_k=N^{(j)}$ for
every $j$. In this way, every municipality and every household order keeps,
up to rounding, the class proportions observed among the heads of household in
the PNS.

\paragraph{Order of generation.} Variables are generated sequentially, so
that each step may use the variables generated before it as predictors:
\begin{enumerate}[label=(\arabic*)]
\item smoking, alcohol consumption and physical exercise, from age, sex
and race
\item self-rated health, using the three variables of step (1);

\item chronic conditions: diabetes, high cholesterol, heart disease,
stroke, asthma, arthritis/rheumatism, depression and
chronic kidney disease, using the variables of steps (1) and (2). Smoking, alcohol
consumption and physical inactivity are established modifiable risk factors
for stroke \cite{odonnell2010}, and physical inactivity is a
well-documented risk factor for depression \cite{schuch2018}. For asthma and
chronic kidney disease, whose association with these three lifestyle
variables is comparatively weak in the literature, only age, sex and race are
used as predictors.

\end{enumerate}
In the validation study (Section~\ref{sec:experiments}), the predictors from earlier steps take
their observed PNS values in both the training and the test sets, so each
variable is evaluated on its own, without propagating errors from earlier
steps. In the synthetic population, only generated values are available, and
the predictors from earlier steps take the values generated in those steps.


\section{Methods}
\label{sec:methods}

\subsection{Notation}

Let $X_\ell=\{(x^i,y^i)\}_{i=1}^{n}\subset\R^d\times\cK$ be the labeled
data (PNS respondents), where $x^i$ is the vector of predictors of
individual $i$ and $y^i$ its observed class, and let
$X_u=\{x^1,\dots,x^m\}\subset\R^d$ be the unlabeled data (for instance, the
synthetic individuals), whose classes are unknown. For a given outcome
with $K$ classes, we write $\cK:=\{1,\dots,K\}$ and $[m]:=\{1,\dots,m\}$.
For each unlabeled point $x^i$, $i\in[m]$, the vector
$p^i=(p^i_1,\dots,p^i_K)$ denotes the estimated probabilities of $x^i$
belonging to each class, with $p^i_k\ge0$ and $\sum_{k\in\cK}p^i_k=1$;
these probabilities are obtained as described in
Sections~\ref{sec:base} and~\ref{sec:ensemble}.
\subsection{Base classifiers}
\label{sec:base}

We use three classifiers, chosen because they rely on different
assumptions about the relation between the predictors $x$ and the class
$y$: multinomial logistic regression assumes linear class scores, the
random forest partitions the predictor space into regions, and the support
vector machine builds nonlinear boundaries through a kernel. Each of them
provides, for every individual, an estimate of $\Pr(y=k\mid x)$,
$k\in\cK$, and these estimates are combined  as shown in
Section~\ref{sec:ensemble}.
\paragraph{Multinomial logistic regression.}
Let $\{(x^i,y^i)\}_{i=1}^n$ be the training data, with $x^i\in\mathbb{R}^d$
and $y^i\in\cK=\{1,\dots,K\}$. With parameters $\theta=(W,b)$,
$W=[w_1,\dots,w_K]$, and linear scores $s_k(x)=w_k^\top x+b_k$, class
probabilities are given by the softmax
\begin{equation}
\Pr(y=k\mid x;\theta)=\frac{\exp(s_k(x))}{\sum_{k'=1}^K\exp(s_{k'}(x))},
\qquad k\in\cK,
\label{eq:softmax}
\end{equation}
and $\theta$ minimizes the regularized cross-entropy
\begin{equation}
J(\theta)=\sum_{i=1}^n\Bigl(-s_{y^i}(x^i)+\log\sum_{k=1}^K e^{s_k(x^i)}\Bigr)
+\frac{\alpha}{2}\lVert W\rVert_F^2,\qquad \alpha>0
\label{eq:loss}
\end{equation}

\cite{hosmer2013}.

Since log-sum-exp is convex and the scores are affine in $\theta$, $J$
is convex, and every local minimizer is global. 
The $\ell_2$ penalty makes $J$ strictly convex in $W$,
which yields a unique minimizer $W^\star$. We use the default of MLJLinearModels.jl, $\alpha=\varepsilon_{\mathrm{mach}}
\approx 2.2\times10^{-16}$, with the intercept not penalized; the penalty then
only guarantees a unique minimizer $W^\star$ and has no practical shrinkage
effect.

\paragraph{Random forest.} A random forest is a collection of $T$ decision
trees, each grown on a random sample, drawn with replacement, of the
training data, with a random subset of predictors considered at each split
\cite{breiman2001}. Given a point $x$, each tree assigns it the majority
class of the leaf it reaches, and the probability of class $k$ is the
fraction of the $T$ trees that assign $x$ to class $k$,
\begin{equation}
p_{k,\mathrm{rf}}(x)=\frac{1}{T}\,\big|\{t\in\{1,\dots,T\}:\ \text{tree } t
\text{ assigns } x \text{ to class } k\}\big|.
\label{eq:rf}
\end{equation}
In particular, $p_{k,\mathrm{rf}}(x)=0$ whenever no tree assigns $x$ to class
$k$, which is frequent for rare classes; this is why the positivity of the
ensemble relies on the logistic regression (Section~\ref{sec:ensemble}).
The hyperparameters are given in Section~\ref{sec:experiments}.

\paragraph{Support vector machine.}
For two classes labeled $y\in\{-1,+1\}$, the support vector machine
\cite{cortes1995} finds the separating hyperplane with the largest margin,
allowing some points to violate it at a cost $C>0$:
\begin{equation}
\min_{w,b,\xi}\ \tfrac12\|w\|^2+C\sum_i \xi_i
\quad\text{s.t.}\quad y^i\big(w^\top\phi(x^i)+b\big)\ge 1-\xi_i,\quad
\xi_i\ge0.
\end{equation}
Here $\phi$ maps the features into a higher-dimensional space, which allows
nonlinear decision boundaries in the original space. In practice, the problem
is solved in its dual form, where $\phi$ appears only through inner products
$\phi(x)^\top\phi(x')$; these are replaced by the RBF kernel
$\kappa(x,x')=\exp(-\gamma\|x-x'\|^2)$, so $\phi$ never needs to be computed
explicitly. The problem is a convex quadratic program, so every local minimum
is global and training does not depend on initialization.

The SVM is binary, so the multiclass case is handled by the one-versus-one
strategy: one classifier is trained for each of the $K(K-1)/2$ pairs of
classes. The SVM also does not output probabilities. For each pair $(j,k)$,
Platt scaling \cite{platt1999} fits a logistic regression on the decision
value, giving an estimate $r_{jk}$ of $\Pr(y=j\mid y\in\{j,k\},x)$. Pairwise
coupling \cite{wu2004} then finds the class probabilities $\Pr(y=k\mid x)$,
$k\in\cK$, summing to one, that are most consistent with all the pairwise
estimates $r_{jk}$.

\subsection{Ensemble}
\label{sec:ensemble}
As shown in Table~\ref{tab:targets}, most of
the variables have strongly imbalanced classes.

In  this setting, classifiers built on different assumptions are known to
behave differently: for instance, SVMs tend to favor the majority class
when classes are imbalanced \cite{akbani2004}, while the behavior of
logistic regression and random forests on rare classes depends on
different aspects of the data. Thus, no single classifier is expected to be
the best one for all variables, and the PNS sample is too small to select a
single model reliably. We therefore combine them: averaging classifiers
whose errors are not identical tends to reduce the variance of the
individual estimates\cite{kuncheva2014}. The combination has a second,
structural advantage: the random forest may assign zero probability to a
rare class (see \eqref{eq:rf}), whereas the logistic regression never
does, so the ensemble probabilities are strictly positive.

Let $p^i_{k,\mathrm{lr}}$, $p^i_{k,\mathrm{rf}}$ and $p^i_{k,\mathrm{svm}}$ be the
probabilities that $x^i$ belongs to class $k$ given by the logistic
regression, the random forest and the SVM, respectively. The ensemble
probability is the average
\begin{equation}
p^i_k=\tfrac13\left(p^i_{k,\mathrm{lr}}+p^i_{k,\mathrm{rf}}+p^i_{k,\mathrm{svm}}\right),
\qquad k\in\cK.
\label{eq:ensemble}
\end{equation}
We use equal weights because no classifier dominates the others across
all outcomes.
Since each term is a probability vector, $\sum_k p^i_k=1$. Moreover, because
$p^i_{k,\mathrm{lr}}>0$ for all $i,k$, the ensemble satisfies
$p^i_k\ge \tfrac13 p^i_{k,\mathrm{lr}}>0$ for all $i,k$, even when the random
forest assigns zero probability to some class. The ensemble prediction (majority vote) is $\hat y^i=\arg\max_k p^i_k$.

\subsection{Classification assignment}
\label{sec:opt}

Let $\lambda_k\in\N$ be the target number of points assigned to class $k$,
with $\sum_{k\in\cK}\lambda_k=m$. The model of \cite{burgard2026fair} allows
these target counts to be prescribed within groups of points; here we do not use groups: the problem is solved separately for each
municipality and each household order, with $m$ the number of synthetic
individuals of that order in the municipality and a single target count
$\lambda_k$ for each class $k$, given by~\eqref{eq:lambda} for the heads of
household and by~\eqref{eq:lambdaj} for the other orders. With assignment variables $z^i_k$ and slacks
$\eta_k,\beta_k$ for over- and under-assignment, the problem, already in
its linear relaxation, reads
\begin{align}
\min_{z,\eta,\beta}\quad & -\frac1m\sum_{i=1}^m\sum_{k\in\cK}p^i_k z^i_k
+\frac B2\sum_{k\in\cK}\left(\eta_k+\beta_k\right)
\label{eq:obj}\\
\text{s.t.}\quad & \sum_{k\in\cK} z^i_k=1, && i\in[m], \label{eq:assign}\\
& \sum_{i=1}^m z^i_k-\eta_k+\beta_k=\lambda_k, && k\in\cK, \label{eq:count}\\
& z^i_k\ge0, && i\in[m],\ k\in\cK, \label{eq:zpos}\\
& \eta_k,\beta_k\ge0, && k\in\cK. \label{eq:slackpos}
\end{align}

The first term of \eqref{eq:obj} favors assigning each point to classes
with high probability, while the second penalizes deviations from the
target counts, weighted by $B>0$. In \cite{burgard2026fair}, $z^i_k$ is
additionally fixed to zero whenever $p^i_k=0$. As shown  in
Section~\ref{sec:ensemble} this never happens for our ensemble, so that
constraint is omitted.
Regarding the constraints of the formulation, Constraint~\eqref{eq:assign} guarantees that each individual point $x^i$ is assigned to exactly one class by enforcing that the assignment variables across all classes sum to one ($\sum_{k \in \mathcal{K}} z_k^i = 1$). In an integral assignment, if an observation belongs to class $3$, for instance, then $z_3^i = 1$ and $z_k^i = 0$ for all other classes $k \neq 3$. 

Next, Constraint~\eqref{eq:count} aligns the total number of points assigned to class $k$ with the prescribed target $\lambda_k$. Merely selecting the class with the largest predicted probability for each individual point is insufficient, as it often violates the global cardinality requirements. To accommodate this, the slack variables $\eta_k$ and $\beta_k$ introduce the necessary flexibility: when a class is under-assigned ($\sum_{i=1}^m z_k^i < \lambda_k$), the deficiency is captured by $\beta_k > 0$ with $\eta_k = 0$; conversely, if a class receives an excess of assignments ($\sum_{i=1}^m z_k^i > \lambda_k$), it leads to $\eta_k > 0$ with $\beta_k = 0$. In the ideal scenario where the assignment exactly meets the target count, both slack variables remain zero ($\eta_k = \beta_k = 0$). 

Finally, Constraints~\eqref{eq:zpos} and~\eqref{eq:slackpos} enforce the non-negativity of both the assignment variables ($z_k^i \ge 0$) and the slack variables ($\eta_k, \beta_k \ge 0$), thereby completing the continuous linear relaxation of the assignment problem.

\begin{theorem}[\cite{burgard2026fair}, Theorem 1]
\label{thm:tu}
The constraint matrix of \eqref{eq:assign}--\eqref{eq:count} is totally
unimodular. Hence every extreme point of the feasible set of
\eqref{eq:obj}--\eqref{eq:slackpos} is integral, and there is an optimal
solution with $z^i_k\in\{0,1\}$ for all $i,k$.
\end{theorem}

\subsection{Structure of the optimal solutions}
\label{sec:structure}

The following result shows that, in our setting, the soft constraints
\eqref{eq:count} are satisfied exactly, and clarifies the role of $B$.
\begin{proposition}
\label{prop:exact}
Assume that $p^i_k>0$ for all $i\in[m]$, $k\in\cK$, that
$\sum_{k\in\cK}\lambda_k=m$, and that $B>1/m$. Then every optimal solution
$(z^*,\eta^*,\beta^*)$ of \eqref{eq:obj}--\eqref{eq:slackpos} satisfies
$\eta^*=\beta^*=0$. Consequently, $z^*$ is an optimal solution of the
transportation problem
\begin{equation}
\max_{z\ge0}\ \sum_{i=1}^m\sum_{k\in\cK}p^i_kz^i_k
\quad\text{s.t.}\quad
\sum_{k\in\cK}z^i_k=1\ (i\in[m]),\qquad
\sum_{i=1}^m z^i_k=\lambda_k\ (k\in\cK),
\label{eq:transport}
\end{equation}
which does not depend on $B$.
\end{proposition}

\begin{proof}
Let $(z,\eta,\beta)$ be optimal. First, $\eta_k\beta_k=0$ for all $k$:
otherwise, decreasing both by $\min\{\eta_k,\beta_k\}>0$ keeps feasibility
and strictly decreases \eqref{eq:obj}, since $B>0$. Let
$c_k:=\sum_{i=1}^m z^i_k$ be the number of points assigned to class $k$.
By \eqref{eq:count} and the previous observation,
$\eta_k=\max\{0,c_k-\lambda_k\}$ and $\beta_k=\max\{0,\lambda_k-c_k\}$.

Suppose, for contradiction, that $c_k\neq\lambda_k$ for some $k$. By
\eqref{eq:assign}, $\sum_k c_k=m=\sum_k\lambda_k$, so there exist classes
$a,b$ with $c_a>\lambda_a$ and $c_b<\lambda_b$. Since $c_a>0$, there is a
point $i$ with $z^i_a>0$. Let
$\varepsilon:=\min\{z^i_a,\ c_a-\lambda_a,\ \lambda_b-c_b\}>0$ and move
$\varepsilon$ of the assignment of point $i$ from class $a$ to class $b$:
$\tilde z^i_a=z^i_a-\varepsilon$, $\tilde z^i_b=z^i_b+\varepsilon$, all
other entries unchanged. This keeps \eqref{eq:assign} and \eqref{eq:zpos},
and \eqref{eq:count} holds with $\tilde\eta_a=\eta_a-\varepsilon\ge0$,
$\tilde\beta_b=\beta_b-\varepsilon\ge0$ and all other slacks unchanged. The
change in the objective is

\[
\left( -\frac{1}{m} \left( p_b^i (z_b^i + \varepsilon) + p_a^i (z_a^i - \varepsilon) \right)
+ \frac{B}{2} \left( (\eta_a - \varepsilon) + (\beta_b - \varepsilon) \right) \right)
- \left( -\frac{1}{m} \left( p_a^i z_a^i + p_b^i z_b^i \right)
+ \frac{B}{2} \left( \eta_a + \beta_b \right) \right)
\]

\[
= -B\varepsilon - \frac{1}{m} \left( p_b^i \varepsilon - p_a^i \varepsilon \right) = -\frac{\varepsilon}{m} \left( p_b^i - p_a^i \right) - B\varepsilon
\]
\[
= -\frac{\varepsilon}{m} \left( p_b^i - p_a^i \right) - B\varepsilon < \frac{\varepsilon}{m} p_a^i - B\varepsilon
\]
\[
\le \varepsilon \left( \frac{1}{m} - B \right)
\]

Since we assumed $B > 1/m$, we have:

\[
\varepsilon \left( \frac{1}{m} - B \right) < 0
\]

That is, the new solution has a smaller objective value, so the original one was not optimal. The assumption $c_k \neq \lambda_k$ is false, therefore $c_k = \lambda_k$ and $\eta_k = \beta_k = 0$ for all $k$.


With zero slacks, \eqref{eq:obj} equals $-\frac1m\sum_{i,k}p^i_kz^i_k$ and
the constraints reduce to those of \eqref{eq:transport}, which proves the
second claim.
\end{proof}

Note that, with $B=2$, as in \cite{burgard2026fair}, the condition $B>1/m$ holds for any $m\ge1$. Note that summing the constraints of \eqref{eq:transport} over $i$ and over $k$ gives $\sum_k\lambda_k=m$; hence, without slack variables, the model is infeasible whenever this identity fails. In \cite{burgard2026fair}, the targets $\lambda_k$ were not required to sum to $m$, so the slacks $\eta,\beta$ were necessary to guarantee feasibility, and the penalty $B$ governed the trade-off between assignment quality and deviation from the targets.

\section{Experimental design}
\label{sec:experiments}

 All experiments are implemented in \texttt{Julia}
1.12.5 \cite{bezanson2017}. The classifiers are fitted with the \texttt{MLJ}
framework \cite{blaom2020mlj}, using \texttt{MLJLinearModels.jl}
\cite{mljlinearmodels} for the logistic regression, \texttt{DecisionTree.jl}
\cite{sadeghi2022} for the random forest  and \texttt{LIBSVM.jl}, a Julia wrapper of LIBSVM
\cite{chang2011}, for the SVM. The linear program \eqref{eq:transport} is modeled in
\texttt{JuMP} \cite{dunning2017} and solved with HiGHS \cite{huangfu2018}.
All experiments were run on a single workstation (Intel Xeon W-2133,
3.6\,GHz, 32\,GB RAM, Ubuntu 20.04). The code is available at \url{https://github.com/vitoriablau/SigeSC}.

\subsection{Validation on PNS data}

The validation mimics the application. A labeled part of the PNS sample for
Santa Catarina plays the role of the survey, and a disjoint held-out part
plays the role of the unlabeled population, whose labels are used \emph{only}
for evaluation. The procedure below is carried out separately for each target
variable in Table~\ref{tab:targets}.
\begin{enumerate}[label=\arabic*.]
\item Split the sample into $5$ folds, stratified by the target variable, and
repeat this with $5$ different seeds ($5\times5$ repeated cross-validation).
Each fold serves once as the test set, with the remaining folds as the
training set, giving $25$ training/test pairs.
\item On each training set, fit the classifiers used by the methods below.
\item On the corresponding test set, let $m$ be its size and set $\lambda_k$
to the number of test units in class $k$, computed from the true labels.
Since every test unit belongs to exactly one class,
$\sum_{k\in\cK}\lambda_k=m$.
\item Apply each of the following methods to the test set and compare their
assignments with the true labels:
\begin{itemize}
\item \textbf{LR:} 
Multinomial logistic regression with an $\ell_2$-regularized multinomial log-likelihood as in \eqref{eq:loss}. Numeric predictors are standardized and categorical predictors one-hot encoded. Models are trained on a random-oversampled version of the training fold (minority classes resampled with replacement up to the majority class size), as the package exposes no native per-class weighting for this formulation.

\item \textbf{RF:} 
Random forest with $T=300$ trees, maximum depth 8, and a minimum of 15 samples per leaf. The choice of $T=300$ ensures adequate variance reduction and convergence of ensemble estimates with negligible additional gain from larger values. The structural constraints (depth 8, minimum leaf size 15) serve as explicit regularizers to prevent overfitting and avoid degenerate (overconfident $0/1$) leaf probability estimates, ensuring well-calibrated class probabilities for downstream aggregation. Trained on the same oversampled fold as LR for consistency.

\item \textbf{SVM:} 
RBF-kernel SVM with cost $C=1$ and bandwidth $\gamma=1/d$ (where $d$ is the feature dimension). The choice $C=1$ provides a balanced trade-off between margin maximization and training penalty, while $\gamma=1/d$ standardizes kernel scale across datasets of varying dimensionality. Trained on the \emph{unresampled} fold with class weights set inversely proportional to class frequencies to preserve natural density geometry while handling imbalance. Multiclass classification is implemented via a one-versus-one decomposition ($K(K-1)/2$ binary classifiers); probabilities are derived using Platt scaling on pairwise decision values followed by pairwise coupling.

\item \textbf{EM:} 
The ensemble method: probability predictions are computed as the unweighted average $p^i_k=\frac13(p^i_{k,\mathrm{lr}}+p^i_{k,\mathrm{rf}}+p^i_{k,\mathrm{svm}})$ as in \eqref{eq:ensemble}, with majority vote prediction $\hat y^i=\arg\max_k p^i_k$. Equal weighting assigns balanced reliance across linear, kernelized, and tree-based decision surfaces.

\item \textbf{EC:} 
Solves problem \eqref{eq:transport} using target counts $\lambda_k$ and class probabilities estimated via the ensemble average in \eqref{eq:ensemble}. Because the LR model yields strictly positive probabilities, the ensemble probability $p^i_k$ remains strictly positive for all classes. Together with $\sum_{k\in\cK}\lambda_k=m$, this guarantees that the hypotheses of Proposition~\ref{prop:exact} hold: slack variables vanish at every optimum of \eqref{eq:obj}--\eqref{eq:slackpos}, enabling the direct solution of \eqref{eq:transport} independently of $B$.
\end{itemize}
\end{enumerate}
\subsection{Metrics}
Let $m$ denote the total number of evaluation units in the test set. For each unit $i \in \{1, \dots, m\}$, let $y^i \in \{1, \dots, K\}$ represent its true class label (ground truth) and $\hat{y}^i \in \{1, \dots, K\}$ represent the predicted class label generated by a given method.

\textit{Accuracy} is the simplest summary (higher is better): the fraction of units the method gets right,
\begin{equation}
\text{Accuracy} \;=\; \frac{1}{m}\sum_{i=1}^{m}\mathbbm{1}\{\hat y^i = y^i\} \in [0,1].
\end{equation}
It weights every unit equally, so each class counts in proportion to its size. This is exactly the problem in imbalanced data like ours: a method that always predicts the majority class can still score a high accuracy.

To see \emph{which} classes drive that number, we look at each class $k$ on its own, treating it as ``positive'' and every other class as ``negative'' (one-vs-rest), and count four quantities:
\begin{equation}
\mathrm{TP}_k=\sum_{i=1}^m \mathbbm{1}\{\hat y^i=k,\,y^i=k\}, \qquad
\mathrm{FP}_k=\sum_{i=1}^m \mathbbm{1}\{\hat y^i=k,\,y^i\neq k\},
\end{equation}
\begin{equation}
\mathrm{FN}_k=\sum_{i=1}^m \mathbbm{1}\{\hat y^i\neq k,\,y^i=k\}, \qquad
\mathrm{TN}_k=\sum_{i=1}^m \mathbbm{1}\{\hat y^i\neq k,\,y^i\neq k\}.
\end{equation}
$\mathrm{TP}_k$ (true positive) is a unit of class $k$ that the method correctly calls $k$; $\mathrm{FP}_k$ (false positive) is a unit the method calls $k$ but that is really something else; $\mathrm{FN}_k$ (false negative) is a unit that really is $k$ but the method missed; $\mathrm{TN}_k$ (true negative) is a unit that is neither truly $k$ nor called $k$. Accuracy above is just $\frac1m\sum_k \mathrm{TP}_k$: every correctly classified unit contributes one true positive to its own class.

From these counts we define, for each class $k$,
\begin{equation}
\text{Precision}_k=\frac{\mathrm{TP}_k}{\mathrm{TP}_k+\mathrm{FP}_k}, \qquad \text{Recall}_k=\frac{\mathrm{TP}_k}{\mathrm{TP}_k+\mathrm{FN}_k}.
\end{equation}
Precision$_k$ answers ``of everything the method labeled $k$, how much really was $k$?''; Recall$_k$ answers ``of everything that really is $k$, how much did the method catch?'' 

$\mathrm{F1}_k$ is the harmonic mean of the two, and macro-F1 is its unweighted average across classes:
\begin{equation}
\mathrm{F1}_k = 2\cdot\frac{\text{Precision}_k\cdot\text{Recall}_k}{\text{Precision}_k+\text{Recall}_k}, \qquad
\text{macro-F1} = \frac1K\sum_{k=1}^{K} \mathrm{F1}_k \in[0,1].
\end{equation}
Taking the unweighted mean over classes, rather than weighting by class size, is what makes macro-F1 (like balanced accuracy) sensitive to performance on rare classes, unlike plain accuracy (for both metrics, higher values indicate better performance).

\medskip\noindent\textbf{Ordinal distance.} Accuracy and macro-F1 both
treat every mistake alike: predicting ``very good'' for someone whose
self-rated health is ``very poor'' costs exactly as much as predicting
``poor''. The three multiclass targets, however, have ordered classes
(smoking: daily, occasionally, not at all; alcohol: never, less than monthly,
monthly or more; self-rated health: very good to very poor), so a prediction
can be wrong by more or by less. We therefore also measure how
far, in number of classes, each prediction lands from the true class,
\begin{equation}
  d^i = \lvert \hat y^i - y^i \rvert, \qquad
  \mathrm{MAE} = \frac{1}{m}\sum_{i=1}^{m} d^i, \qquad
  \mathrm{nMAE} = \frac{\mathrm{MAE}}{K-1} \in [0,1].
  \label{eq:ordinal}
\end{equation}
$d^i=0$ is a correct prediction, and $d^i=K-1$ is the worst possible
mistake. MAE is expressed in number of classes, whereas $\mathrm{nMAE}$ allows comparing targets with different numbers of classes (for both metrics, lower values indicate better performance). For the nine binary targets (physical exercise and the eight chronic conditions), $\mathrm{nMAE}$ reduces to the classification error rate, $1-\mathrm{Accuracy}$; it is reported for them only for uniformity.
Besides the average, we report how many test units fall at each distance
$d=0,1,\dots,K-1$, since two methods with the same MAE can differ in how
often they make large mistakes.

\section{Results}
\label{sec:results}

The experiments have two stages. In the first (Section~\ref{sec:results-pns}), we
use only the PNS sample of Santa Catarina, whose classes are known, to compare
the five methods. For each target, LR, RF and SVM are trained and each produces
a vector of class probabilities; EM assigns the class with the largest average
probability (Eq.~\eqref{eq:ensemble}), and EC reclassifies the same averaged
probabilities by solving problem~\eqref{eq:transport}. In the second stage
(Section~\ref{sec:generation}), we generate the targets for the synthetic
population of Santa Catarina. There the true classes are unknown, so only EC,
the method with the best balance across classes in the first stage, is used:
LR, RF and SVM serve only to produce the probabilities that enter the linear
program, and every individual's class is assigned by EC.

\subsection{Results on PNS data}
\label{sec:results-pns}

Table~\ref{tab:val} reports the mean
and standard deviation of each metric over the 25 training/test pairs
of the $5\times5$ repeated cross-validation on the $N=3{,}738$ PNS respondents aged 15 or older. In each pair, LR and RF are
trained on the oversampled training fold and SVM on the original training fold
with class weights; EC then solves problem~\eqref{eq:transport} on the test fold,
so $m=|X_{\mathrm{test}}|$, about 748 respondents for most targets. For
diabetes and high cholesterol, $N$ is smaller (3{,}539 and 3{,}514) because of
missing or inapplicable answers, and $m$ is about 708 and 703. In this stage, $\lambda_k$ is the number of test
respondents in class $k$. This is an oracle choice: it measures how well EC
assigns the classes when the class counts are right, not how well they are
estimated.

\begin{table*}[ht]
\centering
\caption{Validation on PNS data: mean accuracy, macro-F1 and nMAE over the 25
test sets of the $5\times5$ cross-validation, for the targets in
Table~\ref{tab:targets}. Higher is better for accuracy and macro-F1, lower is
better for nMAE; best value per target and metric in bold. For $K=2$,
$\text{nMAE}=1-\text{Accuracy}$.}
\label{tab:val}
\scriptsize
\setlength{\tabcolsep}{3pt}
\begin{tabular}{lc ccccc ccccc ccccc}
\toprule
& & \multicolumn{5}{c}{Accuracy} & \multicolumn{5}{c}{Macro-F1} & \multicolumn{5}{c}{nMAE} \\
\cmidrule(lr){3-7}\cmidrule(lr){8-12}\cmidrule(lr){13-17}
Target & $K$ & LR & RF & SVM & EM & EC & LR & RF & SVM & EM & EC & LR & RF & SVM & EM & EC \\
\midrule
Smoking & 3 & .362 & .432 & \textbf{.864} & .672 & .777 & .254 & .281 & .309 & .338 & \textbf{.353} & .472 & .433 & \textbf{.130} & .251 & .212 \\
Alcohol & 3 & .523 & .491 & \textbf{.575} & .534 & .538 & .446 & .436 & .406 & \textbf{.450} & .443 & .357 & .363 & .358 & .353 & \textbf{.350} \\
Physical exercise & 2 & .573 & .567 & \textbf{.589} & .586 & .582 & .570 & .562 & .567 & \textbf{.577} & \textbf{.577} & .427 & .433 & \textbf{.411} & .414 & .418 \\
Self-rated health & 5 & .262 & .257 & \textbf{.536} & .342 & .455 & .215 & .207 & .185 & .249 & \textbf{.285} & .309 & .303 & \textbf{.128} & .253 & .162 \\
\midrule
Diabetes & 2 & .671 & .672 & \textbf{.912} & .769 & .863 & .532 & .527 & .477 & .573 & \textbf{.575} & .329 & .328 & \textbf{.088} & .231 & .137 \\
High cholesterol & 2 & .619 & .601 & \textbf{.839} & .655 & .765 & .544 & .537 & .456 & .558 & \textbf{.565} & .381 & .399 & \textbf{.161} & .345 & .235 \\
Heart disease & 2 & .687 & .693 & \textbf{.929} & .794 & .896 & .527 & .522 & .482 & .572 & \textbf{.607} & .313 & .307 & \textbf{.071} & .206 & .104 \\
Stroke & 2 & .698 & .787 & \textbf{.979} & .866 & .962 & .454 & .491 & .495 & .528 & \textbf{.531} & .302 & .213 & \textbf{.021} & .134 & .038 \\
Asthma & 2 & .551 & .674 & \textbf{.943} & .903 & .896 & .404 & .458 & .485 & .512 & \textbf{.515} & .449 & .326 & \textbf{.057} & .097 & .104 \\
Arthritis/rheumatism & 2 & .716 & .699 & \textbf{.921} & .772 & .880 & .561 & .548 & .479 & .584 & \textbf{.591} & .284 & .301 & \textbf{.079} & .228 & .120 \\
Depression & 2 & .587 & .621 & \textbf{.866} & .703 & .796 & .508 & .521 & .464 & .555 & \textbf{.561} & .413 & .379 & \textbf{.134} & .297 & .204 \\
Chronic kidney disease & 2 & .574 & .657 & \textbf{.979} & .900 & .958 & .386 & .415 & .495 & .490 & \textbf{.496} & .426 & .343 & \textbf{.021} & .100 & .042 \\
\bottomrule
\end{tabular}
\end{table*}

SVM has the highest accuracy on every target, but only because
it predicts the majority class. On smoking and on all eight diseases, its macro-F1 equals the
value obtained by always predicting the majority class, and its
accuracy mirrors the majority prevalence (0.979 for stroke, whose
prevalence is 2.1\%). It therefore fails to identify the people who
have the condition, which is what matters for reproducing the
population distribution.

EC is the only method that performs well on all metrics. It has the
highest macro-F1 on 10 of the 12 targets. On physical exercise it ties with EM (0.577),
and on alcohol EM is ahead by 0.007. EC is also second only to SVM in accuracy
on every imbalanced target except asthma. Calibration drives this gain: for
heart disease, it raises accuracy from 0.794 (EM) to 0.896 and macro-F1
from 0.572 to 0.607. Since EM and EC share the same ensemble
probabilities, this implies that incorporating known population targets
is what allows EC to recover minority classes without sacrificing
overall accuracy.

\section{Generation for the synthetic population}
\label{sec:generation}

We generate the target variables for the synthetic population of Santa
Catarina of \cite{pinheiro2026synthpop} with EC, the method that gave
the best balance across classes in Section~\ref{sec:results-pns}.

\paragraph{Heads of household.} The classifiers for the heads of household
(order~$1$) are trained on the $2\,669$ PNS respondents of Santa Catarina who
are heads of household, and the class proportions $\hat\pi_k$ in
\eqref{eq:lambda} are estimated from the same respondents, without sampling
weights.

\paragraph{Other household members.} For household members of order
$j\ge2$, the variables are generated with the same sequential procedure
(Section~\ref{sec:procedure}). To preserve intra-household dependence, the
model for a variable $Y$ uses, as an additional predictor, the value
$Y^{(1)}$ already generated for the head of the household. This model is
trained on the $1\,069$ PNS respondents who are not heads of household. Since
the PNS interviews only one resident per household, the head's class in these
households is not observed; we impute it by applying the order-$1$ model and
EC to the heads, using their age, sex and race/skin color from the household
roster. Only household members aged 15 or older are classified, since the PNS
questionnaire covers this age group.

\paragraph{Assignment and rounding.} EC is solved separately for each
household order $j$ in each of the $295$ municipalities, with $m=N^{(j)}$ and
target counts given by \eqref{eq:lambda} for $j=1$ and by \eqref{eq:lambdaj}
for $j\ge2$. Let $r_k$ denote the unrounded value inside
$\operatorname{round}(\cdot)$, that is, $r_k=\hat\pi_k N^{(1)}$ for $j=1$ and
$r_k=(N^{(j)}/N^{(1)})\lambda^{(1)}_k$ for $j\ge2$; in both cases
$\sum_{k\in\cK}r_k=m$. Since the $r_k$ are not integers in general, we round
them with the largest-remainder method: each $\lambda_k$ starts at
$\lfloor r_k\rfloor$, and the remaining $m-\sum_{k\in\cK}\lfloor r_k\rfloor$
units go to the classes with the largest fractional parts, with ties broken by
the smallest index $k$. The resulting targets differ from $r_k$ by less than
one unit and sum to $m$. Hence every order has, up to rounding, the same class
proportions as the heads of household of its municipality, and the predictor
$Y^{(1)}$ changes which members fall in each class, not how many. Since the
ensemble probabilities are strictly positive and $\sum_k\lambda_k=m$,
Proposition~\ref{prop:exact} allows us to solve \eqref{eq:transport}
directly, without slack variables.

By construction, the class proportions in each municipality equal the
targets. The generated population cannot be validated directly, since
the true classes of its units are unknown.

\section{Conclusion}
\label{sec:conclusion}

We added twelve health variables from the PNS 2019 (smoking, alcohol
consumption, physical exercise, self-rated health and eight chronic conditions)
to the synthetic population of Santa Catarina~\cite{pinheiro2026synthpop}. Class probabilities were
estimated by the average of a multinomial logistic regression, a random forest
and an SVM, and individuals were assigned to classes by the linear program of
Burgard et al.~\cite{burgard2026fair}. We showed that, because the ensemble
probabilities are strictly positive and the target counts sum to the number of
individuals, every optimal solution reproduces the target counts exactly, so
the program reduces to a transportation problem that does not depend on the
penalty parameter. In a $5\times5$ cross-validation on PNS data with the
target counts taken from the held-out labels, the calibrated
classifier (EC) had the highest macro-F1 on 10 of the 12 outcomes, tied with the
ensemble on physical exercise and was 0.007 behind it on alcohol consumption.
The SVM had the highest accuracy on every outcome only because it predicts the
majority class; EC, using the same probabilities as the ensemble, raised
accuracy on every imbalanced outcome except asthma, e.g.\ from 0.794 to 0.896
for heart disease, while recovering the minority classes. In the synthetic
population, the procedure classified all heads of household and all other
members aged 15 or older, solving one linear program per variable for each
municipality and household order; all were solved to optimality and
reproduced the target counts exactly.

Several extensions remain. First, the target counts were derived from
state-level PNS proportions, so all municipalities share the same class
proportions; small-area estimates, or the use of the survey's sampling weights,
would allow regional variation. Second, household members of order $j\ge2$
inherit the class proportions of the heads of household, although they differ
in age composition; proportions estimated by household position or age group
would be more realistic. Third, intra-household dependence enters only through
the head's value of the same variable, and chronic conditions are generated
independently of one another given the lifestyle variables; household-level
predictors and models of the joint distribution of several health variables
could capture these dependencies. Finally, children under 15 are not covered
by the PNS individual questionnaire, and the generated population could be
compared with external aggregates, such as municipal indicators from other
surveys, as an indirect validation.

\section*{Acknowledgments}
We thank the Brazilian agency Conselho Nacional de Desenvolvimento Científico e Tecnológico - CNPq for financial support Grant 444264/2024-8,  314788/2025-5 and  198186/2025-8.

\section*{Use of AI}

During the preparation of this work, the authors used Anthropic's Claude (Claude Opus 5.5)~\cite{anthropic2026claude}
to assist with mathematical and editorial review, consistency checks
between the manuscript and the accompanying code, drafting of text, and
implementation of the models. All AI-assisted text and code were reviewed,
edited, and tested by the authors, who take full responsibility for the
content of this publication.

\end{document}